\documentclass[12pt,twoside]{article}
\usepackage[english]{babel}
\usepackage[latin1]{inputenc}
\usepackage{amsmath}
\usepackage{amssymb,amsfonts}
\usepackage{graphicx}  
\usepackage[colorlinks=true, linkcolor=blue, citecolor=blue, urlcolor=blue]{hyperref} 

\newcommand{\bR}{\mathbf{R}}

\newcommand{\bS}{\mathbf{S}}

\newcommand{\Bb}{\boldsymbol{b}}

\newcommand{\HYP}{\mathbb{H}^3}
\newcommand{\HYN}{\mathbb{H}^n}
\newcommand{\SXR}{\bS^2\!\times\!\bR}

\begin{document}
\pagestyle{myheadings}
\markboth{\centerline{Arnasli Yahya and Jen\H o Szirmai}}
{Decomposition method and upper bound density...}
\title
{Decomposition method and upper bound density related to congruent 
saturated hyperball packings in hyperbolic $n-$space}

\author{Arnasli Yahya and Jen\H o Szirmai \\
\normalsize Department of Algebra and Geometry, Institute of Mathematics,\\
\normalsize Budapest University of Technology and Economics, \\
\normalsize M\H uegyetem rkp. 3., H-1111 Budapest, Hungary \\
\normalsize arnasli@math.bme.hu,~szirmai@math.bme.hu
\date{\normalsize{\today}}}

\maketitle


\begin{abstract}
In this paper, we study the problem of hyperball (hypersphere) packings in 
$n$-dimensional hyperbolic space ($n \ge 4$). We prove that to each $n$-dimensional congruent saturated hyperball packing, there is 
an algorithm to obtain a decomposition of $n$-dimensional hyperbolic space $\HYN$ into truncated
simplices. We prove, using the above method and the results of the paper \cite{M94}, that the upper bound of the density for saturated congruent hyperball packings, related to the corresponding truncated tetrahedron cells, is attained in a regular truncated simplex. In 4-dimensional hyperbolic space, we determined this upper bound density to be approximately $0.75864$.
Moreover, we deny A.~Przeworski's conjecture \cite{P13} regarding the monotonization of the density function in the $4$-dimensional hyperbolic space.
\end{abstract}

\newtheorem{theorem}{Theorem}[section]
\newtheorem{cor}[theorem]{Corollary}
\newtheorem{conjecture}{Conjecture}[section]
\newtheorem{lemma}[theorem]{Lemma}
\newtheorem{exmple}[theorem]{Example}
\newtheorem{defn}[theorem]{Definition}
\newtheorem{rmrk}[theorem]{Remark}
\newenvironment{definition}{\begin{defn}\normalfont}{\end{defn}}
\newenvironment{remark}{\begin{rmrk}\normalfont}{\end{rmrk}}
\newenvironment{example}{\begin{exmple}\normalfont}{\end{exmple}}
\newenvironment{acknowledgement}{Acknowledgement}


\section{Introduction, preliminarly results}
In $n$-dimensional hyperbolic geometry occur several new questions 
concerning the  packing and covering problems, e.g.,
in $\HYN$ there are $3$ types of "generalized balls (spheres)":  
the {\it usual balls} (spheres),
{\it horoballs} (horospheres) and {\it hyperballs} (hyperspheres). 
Moreover, the definition of packing or covering density is crucial in hyperbolic spaces
as shown by B\"or\"oczky \cite{B78}, for standard examples see \cite{G--K--K}, \cite{Fe23} and \cite{R06}.
The most widely accepted notion of packing density considers the local densities of balls with respect to their Dirichlet-Voronoi cells
(cf. \cite{B78} and \cite{K98}). In order to consider ball packings in 
$\overline{\mathbb{H}}^n$, we use the extended notion of such local density.

In space $X^n$, let $d_n(r) ~(n\geq 2)$ be the density of $n+1$ mutually touching spheres 
or horospheres of radius $r$ (in case of horosphere $r=\infty$) with respect
to the simplex spanned by their centres. L.~Fejes T\'oth and H.~S.~M.~Coxeter
conjectured that the packing density of balls of radius $r$ in $X^n$ cannot exceed $d_n(r)$.
This conjecture has been proved by C.~A.~Rogers for Euclidean space $\mathbb{E}^n$.
The 2-dimensional spherical case was settled by L.Fejes T\'oth in \cite{FTL}. 
The greatest density in $\mathbb{H}^3$ is $\approx 0.85328$
which is not realized by packing with equal balls. However, 
it is attained by the horoball packing
(in this case $r=\infty$) of
$\overline{\mathbb{H}}^3$ where the ideal centers of horoballs lie on the
absolute figure of $\overline{\mathbb{H}}^3$. This ideal regular
tetrahedron tiling is given with Coxeter-Schl\"afli symbol $\{3,3,6\}$.
Ball packings of hyperbolic $n$-space and of other Thurston geometries
are extensively discussed in the literature see e.g., \cite{Be}, \cite{B78},
\cite{G--K--K}, \cite{MSz17} and \cite{Sz14-1}, where the reader could find further references as well. However, there are countless open questions between ball and horoball packings and coverings, and their investigation is ongoing
(see \cite{KSz}, \cite{KSz14}, \cite{KSz23}, \cite{MSz17}, \cite{MSz18}, \cite{Boh}, \cite{MSz24}, \cite{YSz22}). 

{\it In this paper, we study the question related to saturated, congruent hyperball packings: 
What are the optimal congruent hyperball packing configurations and how large is their density ($ n\ge 4$)?}

In the hyperbolic plane $\mathbb{H}^2$, the universal upper bound of the hypercycle packing density is $\frac{3}{\pi}$,
proved by I.~Vermes in \cite{V79}, and the universal lower bound of the hypercycle covering density is $\frac{\sqrt{12}}{\pi}$,
determined by I.~Vermes in \cite{V81}. 

In \cite{Sz06-1} and \cite{Sz06-2}, we have studied the regular prism tilings (simple truncated Coxeter orthoscheme tilings) and the corresponding optimal hyperball packings in
$\mathbb{H}^n$ $(n=3,4)$. Moreover, we extended the method developed in the former paper \cite{Sz13-3} to 
5-dimensional hyperbolic space, and constructed to each investigated Coxeter tiling a regular prism tiling.
We have studied the corresponding optimal hyperball packings by congruent hyperballs.
Moreover, their metric data and their densities have been determined. In paper \cite{Sz13-4}, we have studied the $n$-dimensional hyperbolic regular prism honeycombs
and the corresponding coverings by congruent hyperballs and we determined their least dense covering densities.
Moreover, we have formulated a conjecture for the candidate of the least dense hyperball
covering by congruent hyperballs in the 3- and 5-dimensional hyperbolic space ($n \in \mathbb{N},3 \le n \le 5)$.

In \cite{Sz17-2}, we discussed congruent and non-congruent hyperball (hypersphere) packings of the truncated regular tetrahedron tilings.
These are derived from the Coxeter simplex tilings $\{p,3,3\}$ $(7\le p \in \mathbb{N})$ and $\{5,3,3,3,3\}$
in $3$- and $5$-dimensional hyperbolic space.
We determined the densest hyperball packing arrangement and its density
with congruent hyperballs in $\mathbb{H}^5$ and determined the smallest density upper bounds of 
non-congruent hyperball packings generated by the above tilings in $\HYN,~ (n=3,5)$.

In \cite{Sz17-1}, we deal with the packings derived by horo- and hyperballs (briefly hyp-hor packings) in $n$-dimensional hyperbolic spaces $\HYN$
($n=2,3$) which form a new class of the classical packing problems.
We constructed in the $2-$ and $3-$dimensional hyperbolic spaces hyp-hor packings that
are generated by complete Coxeter tilings of degree $1$, i.e., the fundamental domains of these tilings are simple frustum orthoschemes. We determined their densest packing configurations and their densities.
We proved, using also numerical approximation methods, that in the hyperbolic plane ($n=2$) the density of the above hyp-hor packings arbitrarily approximate
the universal upper bound of the hypercycle or horocycle packing density $\frac{3}{\pi}$ and
in $\HYP$ the optimal configuration belongs to the $\{7,3,6\}$ Coxeter tiling with density $\approx 0.83267$.
Furthermore, we analyzed the hyp-hor packings in
truncated orthosche\-mes $\{p,3,6\}$ $(6< p < 7, ~ p\in \mathbb{R})$ whose
density function is attained its maximum for a parameter which lies in the interval $[6.05,6.06]$
and the densities for parameters lying in this interval are larger that $\approx 0.85397$. That means that these
locally optimal hyp-hor configurations provide larger densities than the B\"or\"oczky-Florian density upper bound
$(\approx 0.85328)$ for ball and horoball packings but these hyp-hor packing configurations can not be extended to the entirety of hyperbolic space $\mathbb{H}^3$. 

In \cite{Sz17}, we studied a large class of hyperball packings in $\HYP$
that can be derived from truncated tetrahedron tilings.
In order to get a density upper bound for the above hyperball packings, it is sufficient 
to determine this density upper bound locally, e.g., in truncated tetrahedra.
Thus we proved that if the truncated tetrahedron is regular, then the density
of the densest packing is $\approx 0.86338$. This is larger than the B\"or\"oczky-Florian density upper bound for balls and horoballs
but our locally optimal hyperball packing configuration cannot be extended to the entirety of
$\mathbb{H}^3$. However, we described a hyperball packing construction, 
by the regular truncated tetrahedron tiling under the extended Coxeter group $\{3, 3, 7\}$ with maximal density $\approx 0.82251$.

Recently, (to the best of author's knowledge) the candidates for the densest hyperball
(hypersphere) packings in the $3,4$ and $5$-dimensional hyperbolic space $\mathbb{H}^n$ are derived by the regular prism
tilings which have been in papers \cite{Sz06-1}, \cite{Sz06-2}, and \cite{Sz13-3}.

{In \cite{Sz17-2}, we considered hyperball packings in 
$3$-dimensional hyperbolic space. We developed a decomposition algorithm that provides a decomposition of $\HYP$ 
into truncated tetrahedra for each saturated hyperball packing. Therefore, in order to get a density upper bound for hyperball packings, it is sufficient to determine
the density upper bound of hyperball packings in truncated simplices.}

In \cite{Sz18}, we studied hyperball packings related to truncated regular 
octahedron and cube tilings derived from the Coxeter simplex tilings 
$\{p,3,4\}$ $(7\le p \in \mathbb{N})$ and $\{p,4,3\}$ $(5\le p \in \mathbb{N})$
in $3$-dimensional hyperbolic space $\HYP$. We determined the densest hyperball packing arrangement and its density
with congruent and non-congruent hyperballs related to the above tilings. 
Moreover, we proved that the locally densest congruent or non-congruent hyperball configuration belongs to the regular truncated cube with density
$\approx 0.86145$. This is larger than the B\"or\"oczky-Florian density upper bound for balls and horoballs.
We described a non-congruent hyperball packing construction derived from the 
regular cube tiling under the extended Coxeter group $\{4, 3, 7\}$ 
with maximal density $\approx 0.84931$.

In \cite{Sz18-1}, we studied congruent and non-congruent hyperball packings generated by doubly truncated Coxeter orthoscheme tilings in the $3$-dimensional hyperbolic space.
We proved that the densest congruent hyperball packing belongs to the Coxeter orthoscheme tiling of parameter $\{7,3,7\}$ with density $\approx 0.81335$.
This density is equal -- by our conjecture -- with the upper bound density of the corresponding non-congruent hyperball arrangements. 

In \cite{Sz23}, we proved that the density upper bound of the saturated 
congruent hyperball (hypersphere) packings related to the
corresponding truncated tetrahedron cells is realized in a regular truncated tetrahedra with density $\approx  0.86338$.
Furthermore, we proved that the density of locally optimal congruent hyperball arrangement 
in regular truncated tetrahedron is not monotonically increasing function of the height (radius) of corresponding optimal hyperball, 
contrary to the ball (sphere) and horoball (horosphere) packings.

{\it In this paper, we continue and extend the study of 
hypersphere packings to the higher-dimensional hyperbolic space $\HYN$ ($n\ge4$). 
We prove that to each $n$-dimensional congruent, saturated hyperball packing there is 
an algorithm to get a decomposition of $n$-dimensional hyperbolic space into truncated
simplices. We prove that the density upper bound of the saturated congruent hyperball packings related to the
corresponding truncated tetrahedron cells is realized in a regular truncated tetrahedra and in $4$-dimensional hyperbolic space we determined this
upper bound density $\approx 0.75864$.
Moreover, we denies 
A.~Przeworski's conjecture \cite{P13} regarding the monotonization of the congruent hyperball packing density function in the $4$-dimensional hyperbolic space.}

\section{The projective model and saturated hyperball packings of $\HYN$}
For $\mathbb{H}^n$ we use the projective model in the Lorentz space $\mathbb{E}^{1,n}$ of signature $(1,n)$,
i.e.,~$\mathbb{E}^{1,n}$ denotes the real vector space $\mathbf{V}^{n+1}$ equipped with the bilinear
form of signature $(1,n)$
$
\langle  \mathbf{x},\mathbf{y} \rangle = -x^0y^0+x^1y^1+ \dots + x^n y^n
$
where the non-zero vectors
$
\mathbf{x}=(x^0,x^1,\dots,x^n)\in\mathbf{V}^{n+1} \ \  \text{and} \ \ \mathbf{y}=(y^0,y^1,\dots,y^n)\in\mathbf{V}^{n+1},
$
are determined up to real factors, for representing points of $\mathcal{P}^n(\mathbb{R})$. Then $\mathbb{H}^n$ can be interpreted
as the interior of the quadric
$
Q=\{[\mathbf{x}]\in\mathcal{P}^n | \langle  \mathbf{x},\mathbf{x} \rangle =0 \}=:\partial \mathbb{H}^n
$
in the real projective space $\mathcal{P}^n(\mathbf{V}^{n+1},
\mbox{\boldmath$V$}\!_{n+1})$.

The points of the boundary $\partial \mathbb{H}^n $ in $\mathcal{P}^n$
are called points at infinity of $\mathbb{H}^n $, the points lying outside $\partial \mathbb{H}^n $
are said to be outer points of $\mathbb{H}^n $ relative to $Q$. Let $P([\mathbf{x}]) \in \mathcal{P}^n$, a point
$[\mathbf{y}] \in \mathcal{P}^n$ is said to be conjugate to $[\mathbf{x}]$ relative to $Q$ if
$\langle  \mathbf{x},\mathbf{y} \rangle =0$ holds. The set of all points which are conjugate to $P([\mathbf{x}])$
form a projective (polar) hyperplane
$
pol(P):=\{[\mathbf{y}]\in\mathcal{P}^n | \langle  \mathbf{x},\mathbf{y} \rangle =0 \}.
$
Thus the quadric $Q$ induces a bijection
(linear polarity $\mathbf{V}^{n+1} \rightarrow
\mbox{\boldmath$V$}\!_{n+1})$)
from the points of $\mathcal{P}^n$
onto its hyperplanes.

The point $X [\bold{x}]$ and the hyperplane $\alpha [\mbox{\boldmath$a$}]$
are called incident if $\bold{x}\mbox{\boldmath$a$}=0$ ($\bold{x} \in \bold{V}^{n+1} \setminus \{\mathbf{0}\}, \ \mbox{\boldmath$a$} \in \mbox{\boldmath$V$}_{n+1}
\setminus \{\mbox{\boldmath$0$}\}$).

The hypersphere (or equidistance surface) is a quadratic surface at a constant distance
from a plane (base plane) in both halfspaces. The infinite body of the hypersphere, containing the base plane, is called hyperball.

The {\it $n$-dimensional half hyperball} with distance $h$ to a base plane $\beta$ is denoted by $\mathcal{H}^{h^+}_n$.
The volume of a bounded hyperball piece $\mathcal{H}^{h^+}_n(\mathcal{A})$,
delimited by a $(n-1)$-polygon $\mathcal{A}_{n-1} \subset \beta$. In the hyperbolic plane its area can be determined by the classical formula 
of J.~Bolyai \cite{B31}. The volume of a bounded hyperball piece $\mathcal{H}^{h^+}_n(\mathcal{A}_{n-1})$
bounded by an $(n-1)$-polytope $\mathcal{A}_{n-1} \subset \beta$, $\mathcal{H}^{h^+}_n$ and by 
hyperplanes orthogonal to $\beta$ derived from the facets of $\mathcal{A}_{n-1}$ can be determined by the formulas (\ref{Eq:2.1}), (\ref{Eq:2.2}) and (\ref{Eq:2.3}) 
that follow from the suitable extension of the classical method of J.~Bolyai:
\begin{equation}\label{Eq:2.1}
Vol_3(\mathcal{H}^{h^+}_3(\mathcal{A}_2))=\frac{1}{4}Vol_2(\mathcal{A}_{2})\left[k \sinh \frac{2h}{k}+
2 h \right], \tag{2.1}
\end{equation}
\begin{equation}\label{Eq:2.2}
Vol_4(\mathcal{H}^{h^+}_4(\mathcal{A}_3))=\frac{1}{8} Vol_3(\mathcal{A}_{3})k \left[ \frac{2}{3} \sinh \frac{3h}{k}+
6 \sinh \frac{h}{k} \right], \tag{2.2}
\end{equation}
\begin{equation}\label{Eq:2.3}
Vol_5(\mathcal{H}^{h^+}_5(\mathcal{A}_4))=\frac{1}{16} Vol_4(\mathcal{A}_4)\left[ k  \left( \frac{1}{2} \sinh \frac{4h}{k}+
4 \sinh \frac{2h}{k}\right) +6 h \right], \tag{2.3}
\end{equation}
where the volume of the hyperbolic $(n-1)$-polytope $\mathcal{A}_{n-1}$ lying in the plane $\beta$ is $Vol_{n-1}(\mathcal{A}_{n-1})$.
The constant $k =\sqrt{\frac{-1}{K}}$ is the natural length unit in
$\mathbb{H}^n$, where $K$ denotes the constant negative sectional curvature. In the following we may assume that $k=1$. 

Let $\mathcal{H}^h$ be a hyperball packing in $\HYN$ with congruent hyperballs of height $h$ in the $n$-dimensional hyperbolic space $\HYN$. 

The notion of {\it saturated packing} follows from that fact that the density of any packing can be improved by adding further packing elements as long as there is 
sufficient room to do so. However, we usually apply this notion for packings with congruent elements.
Now, we modify the classical definition of saturated packing for non-compact ball packings in $n$-dimensional hyperbolic space $\HYN$ ($n\ge 2$ integer parameter) (see \cite{Sz17-2}):
\begin{defn}\label{def:2.1}
A ball packing with non-compact balls (horoballs or/and hyperballs) in $\HYN$ is saturated if no new non-compact ball can be added to it. 
\end{defn}
\section{Decomposition into truncated tetrahedra}
We take the set of hyperballs $\{ \mathcal{H}^{h}_i\}$ of a congruent saturated hyperball packing $\mathcal{H}^h$ (see Definition \ref{def:2.1}) in the $n$-dimensional hyperbolic space $\HYN$ ($3<n \in \mathbb{N}$).
Their base hyperplanes are denoted by $\beta_i$.
Thus in a saturated hyperball packing the distance between two ultraparallel base hyperplanes
$d(\beta_i,\beta_j)$ is at least $2h$ (where for the natural indices holds $i < j$
and $d$ is the hyperbolic distance function).

In this section we describe a procedure to get a decomposition of $n$-dimensional hyperbolic space $\HYN$ ($3<n \in \mathbb{N}$) into truncated
$n$-simplices corresponding to a given saturated hyperball packing.
\begin{enumerate}
\item The notion of the radical hyperplane (or power hyperplane) of two Euclidean spheres can be extended to the hyperspheres.
The radical hyperplane (or power hyperplane) of two non-intersecting hyperspheres
is the locus of points at which tangents drawn to both hyperspheres have the same length (so these points have equal power with respect
to the two non-intersecting hyperspheres). If the two non-intersecting hyperspheres are congruent also in Euclidean sense in the model then their 
radical hyperplane coincides with their "Euclidean symmetry hyperplane" and any two congruent hypersphere can be transformed into such an hypersphere arrangement.

Using the radical hyperplanes to the hyperballs $\mathcal{H}^h_i$ similarly to the Euclidean space, can be constructed
the unique Dirichlet-Voronoi (in short $D-V$) decomposition of $\HYN$ to the given hyperball packing $\mathcal{B}^h$. Now, the
$D-V$ cells are infinite hyperbolic polytopes containing the corresponding hyperball,
and its vertices are proper points of $\mathbb{H}^n$. We note here (it is easy to see), that a vertex of any $D-V$ cell 
cannot be outer or boundary point of $\mathbb{H}^n$ relative to $Q$, because the hyperball packing $\mathcal{B}^h$ is saturated by the Definition \ref{def:2.1}. 

\item We consider an arbitrary {\it proper} vertex $P \in \HYN$ of the above $D-V$ decomposition and the hyperballs $\mathcal{H}^h_i(P)$
whose $D-V$ cells
meet at $P$. The base hyperplanes of the hyperballs $\mathcal{H}^h_i(P)$ are denoted by $\beta_i(P)$, and these hyperplanes determine a
non-compact polytop $\mathcal{D}^i(P)$ with the intersection of their halfspaces 
containing the vertex $P$. Moreover, denote $A_1,A_2,A_3,\dots$ the outer vertices of $\mathcal{D}^i(P)$ and cut off 
$\mathcal{D}^i(P)$ with the polar hyperplanes $\alpha_j(P)$ of its outer vertices $A_j$. Thus, we obtain a convex compact polytop 
$\mathcal{D}(P)$.
This is bounded by the base planes $\beta_i(P)$ and "polar hyperplanes" $\alpha_i(P)$. Applying this procedure for all vertices of the
above Dirichlet-Voronoi decomposition, we obtain an other decomposition of $\HYN$ into convex polytopes.
\item We consider $\mathcal{D}(P)$ as a tile of the above decomposition. The hyperplanes from the finite set of base hyperplanes $\{\beta_i(P)\}$ are 
called adjacent if there is a vertex $A_s$ of $\mathcal{D}^i(P)$ that lies on each of the above hyperplanes.
We consider non-adjacent planes $\beta_{k_1}(P),\beta_{k_{2}}(P),\beta_{k_{3}}(P), \dots \beta_{k_{m}}(P) \in \{\beta_i(P)\}$ $(k_l  \in \mathbb{N}^+, ~ l=1,2,3,\dots m)$
that have an outer point of intersection denoted by $A_{k_1\dots k_m}$. Let $N_{\mathcal{D}(P)}$ denote the {\it finite} number of the outer points $A_{k_1\dots k_m}$ 
related to $\mathcal{D}(P)$, it is clear, that its 
minimum is $0$ if $\mathcal{D}^i(P)$ is tetrahedron. 
The polar hyperplane $\alpha_{k_1\dots k_m}$ of $A_{k_1\dots k_m}$ is orthogonal to hyperplanes $\beta_{k_1}(P),\beta_{k_2}(P), \dots \beta_{k_m}(P)$ 
(thus it contain the poles $B_{k_1}$, $B_{k_2}$,\dots $B_{k_m}$ of planes $\beta_{k_1}(P),\beta_{k_2}(P),\dots \beta_{k_m}(P)$) and divides $\mathcal{D}(P)$ into two convex polytopes
$\mathcal{D}_1(P)$ and $\mathcal{D}_2(P)$. 
\item If $N_{\mathcal{D}_1(P)} \ne 0$ and $N_{\mathcal{D}_2(P)} \ne 0$ then $N_{\mathcal{D}_1(P)} < N_{\mathcal{D}(P)}$ and $N_{\mathcal{D}_2(P)} < N_{\mathcal{D}(P)}$. 
We apply the point 3 for polytopes $\mathcal{D}_i(P),~ i \in \{1,2\}$.
\item If $N_{\mathcal{D}_i(P)} \ne 0$ or $N_{\mathcal{D}_j(P)} = 0$ ($i \ne j,~ i,j\in\{1,2\}$) then we consider the polytope $\mathcal{D}_i(P)$ where 
$N_{\mathcal{D}_i(P)}=N_{\mathcal{D}(P)}-1$ because the vertex $A_{k_1\dots k_m}$ is left out and apply the point 3.
\item If $N_{\mathcal{D}_1(P)}=0$ and $N_{\mathcal{D}_2(P)}=0$ then the procedure is over for $\mathcal{D}(P)$. We continue the procedure with the next cell. 
\item It is clear, that the above hyperplane $\alpha_{k_1\dots k_m}$ intersects every hyperball $ \mathcal{H}^h_j(P)$, $(j=k_1\dots k_m)$.

\begin{lemma}
The hyperplane $\alpha_{k_1 \dots k_m}$ of $A_{k_1 \dots k_m}$ does not intersect the hyperballs $\mathcal{H}^h_s(P)$ where $A_{k_1 \dots k_m} \notin \beta_s(P).$
\end{lemma}
{\bf Proof}

Let $\mathcal{H}^h_s(P)$, ($A_{k_1 \dots k_m} \notin \beta_s(P)$) be an arbitrary hyperball corresponding to $\mathcal{D}(P)$ with base plane $\beta_s(P)$ whose pole
is denoted by $B_s$.
The common perpendicular $\sigma$ of the hyperplanes $\alpha_{k_1 \dots k_m}$ and $\beta_s(P)$ is the line through the point $A_{k_1 \dots k_m}$ and $B_s$.
We take a hyperplane $\kappa$ containing the above common perpendicular, and its intersections with $\mathcal{D}(P)$ and $\mathcal{H}^h_s(P)$ are denoted by
$\phi$ and $\eta$. Using the correspondig results of \cite{V79}, \cite{Sz17-2} and applying the complete induction principle we obtain that 
the $(n-1)$-dimensional subspace (hyperplane) $\phi=\alpha_{k_1 \dots k_m}\cap \kappa$ does not intersect the hypershere $\eta=\mathcal{H}^h_s(P)\cap \kappa$.
The hyperplane $\alpha_{k_1 \dots k_m}$ and the hypersphere $\mathcal{H}^h_s(P)$ can be generated by rotation of $\phi$ and $\eta$ about
the common perpendicular $\sigma$; therefore, they are disjoint. $\square$
\item We have seen in the steps 3, 4, and 5 that the number of the vertices $A_{k_1 \dots k_m}$ of any polytop obtained after the cutting process is equal or less than the original one, and
we have proven in step 7 that the original hyperballs form packings in the new polyhedra $\mathcal{D}_1(P)$ and $\mathcal{D}_2(P)$, as well.
We continue the cutting procedure described in step 3 for both polytopes $\mathcal{D}_1(P)$ and $\mathcal{D}_2(P)$. If a derived polyhedron is a truncated
tetrahedron then the cutting procedure does not give new polyhedra, thus the procedure will not be continued.
Finally, after a finite number of cuttings we get a decomposition of $\mathcal{D}(P)$ into truncated simplices, 
and in any truncated simplex the corresponding congruent hyperballs from $\{ \mathcal{H}^h_i\}$ form a packing. Moreover, we apply the above method for the further cells.
\end{enumerate}
Finally we get the following 
\begin{theorem}
The above described algorithm provides for each congruent saturated hyperball packing a decomposition of $\HYN$ into truncated simplices. ~ ~$\square$
\end{theorem}
\begin{rmrk}
In \cite{P13}, Przeworski proved a similar theorem but it was true only for cases if the base planes of hyperspheres form "symmetric cocompacts arrangements" in $\HYN$.
\end{rmrk}
In \cite{M94}, Miyamoto proved the analogue theorem of K.~B\"or\"oczky's theorem in \cite{B78} and \cite{B--F64}:
\begin{theorem}[Y.~Miyamoto, \cite{M94}]\label{thm:Miyamoto}
If a region in $\HYN$ bounded by hyperplanes has a hyperball (hypersphere) packing of height (radius) $r$ about its boundary, then in some sense, 
the ratio of its volume to the volume of its boundary is at least that of a regular truncated simplex of (inner) edgelength $2r$.
\end{theorem}
\begin{rmrk}
Independently from the above paper, Przeworski proved a similar theorem with other methods in \cite{P13}.
\end{rmrk}
From the above section it follows that, to each congruent saturated hyperball packing $\mathcal{B}^h$ of hyperballs $\mathcal{H}^h_i$ there is a
decomposition of $\HYN$ into truncated simplices, moreover, applying the Theorem \ref{thm:Miyamoto} of Miyamoto we obtain the following

\begin{lemma}\label{lemma:3.6} In order to get a density upper bound for the congruent saturated hyperball packing in $\HYN$ $(n \ge 3)$, it is sufficient to determine
the density upper bound of hyperball packings in truncated regular simplices.
\end{lemma}
Based on the previous lemma, we can theoretically give a density upper bound to each congruent saturated hyperball packing in $\HYN$ for all $n\ge4$
(for $n=3$ see \cite{Sz23}).
We determine it in the next section in $\mathbb{H}^4$.
\section{An upper bound of the density for packing of congruent hyperballs in hyperbolic $4-$space}
\begin{figure}[ht]
\centering
\includegraphics[width=5.5cm]{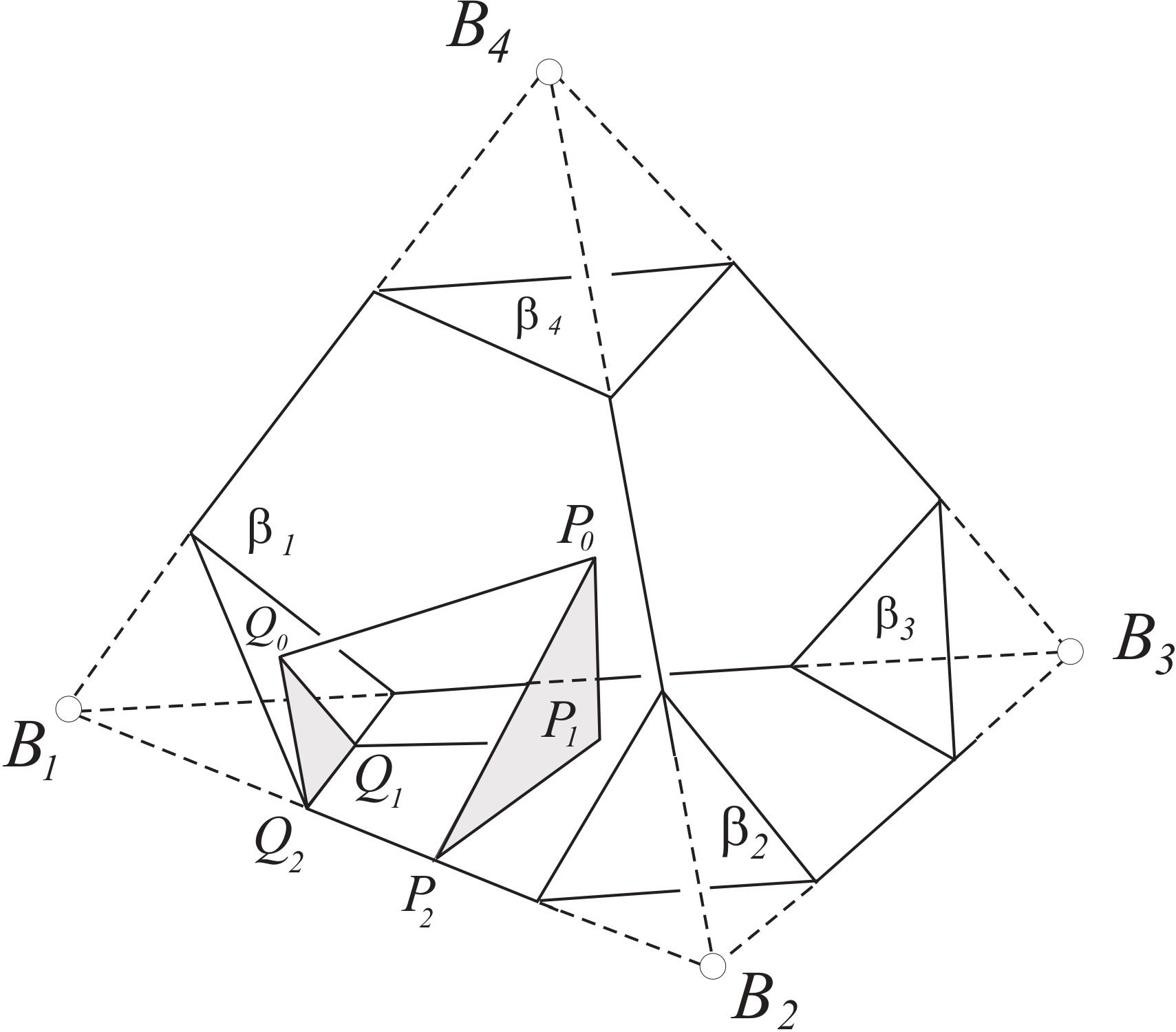} \includegraphics[width=6.5cm]{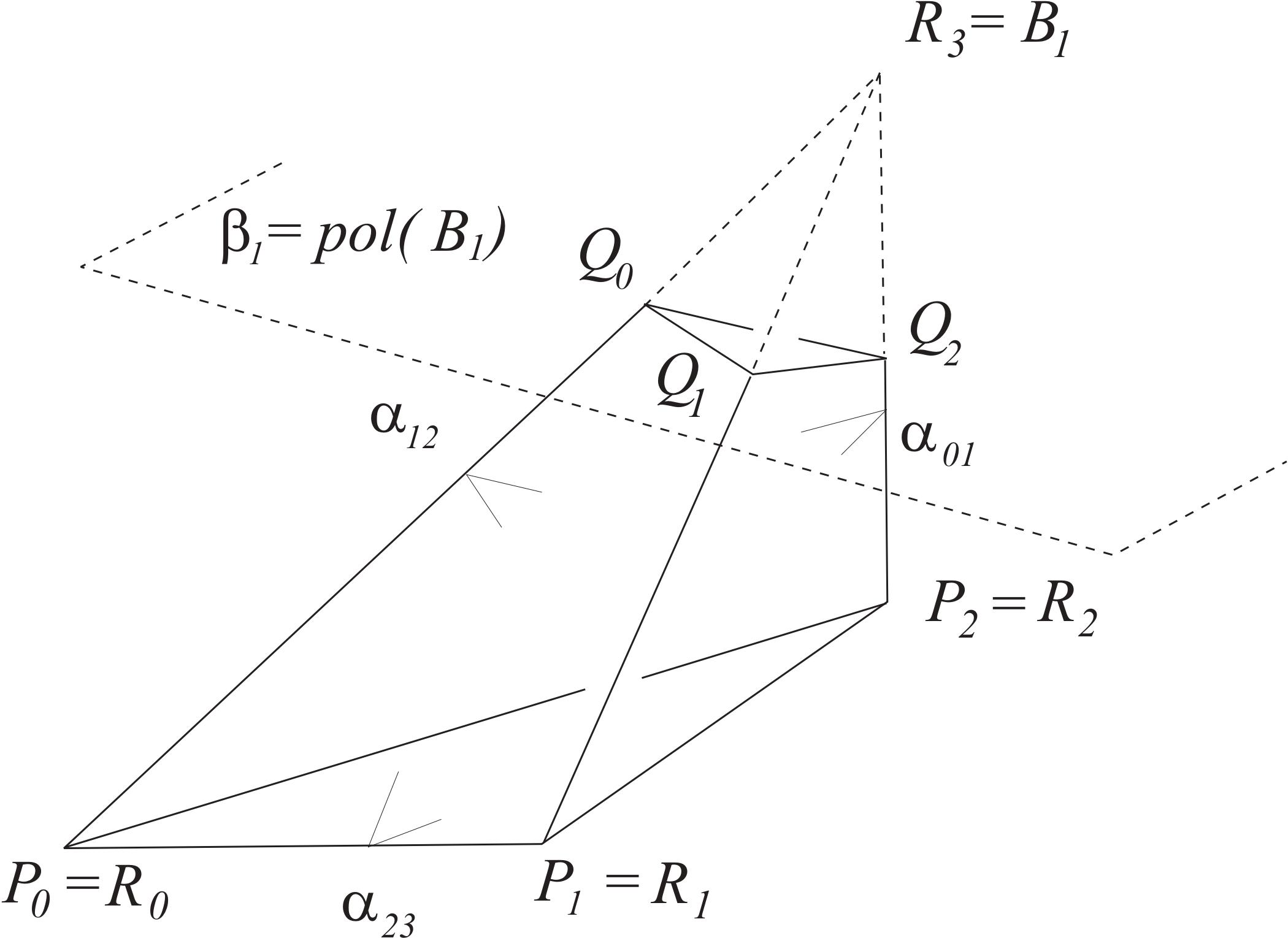}

a.~~~~~~~~~~~~~~~~~~~~~~~~~~~~~~~~~~~~~b.
\caption{Truncated tetrahedron with a complete orthoscheme of degree $m=1$ (simple frustum orthoscheme) in $3$ dimensional hyperbolic space $\HYP$.}
\label{Fig.1}
\end{figure}

From the above section, it follows that to each congruent saturated hyperball packing $\mathcal{B}^h$ of hyperballs $\mathcal{H}^h_i$ in $\HYN$,  there is a
decomposition of $\HYP$ into truncated tetrahedra. Moreover, one of them $\mathcal{S}$, is illustrated in Fig.\ref{Fig.1}.a using the former denotations in the $3$-dimensional hyperbolic space. 

The ultraparallel base planes of $\mathcal{H}^h_i$ $(i=1,2,3,4)$ are denoted by $\beta_i$. The distance between two base planes
$d(\beta_i,\beta_j)=:e_{ij}$ ($i < j$, $i,j \in \{1,2,3,4\})$ and $d$ is the hyperbolic distance function) at least $2h$.
Moreover, the volume of the truncated simplex $\mathcal{S}$ is denoted by $Vol(\mathcal{S})$. We introduce the locally density function
$\delta(\mathcal{S}(h))$ related to $\mathcal{S}$:
\begin{definition}
\begin{equation}
\delta(\mathcal{S}(h)):=\frac{\sum_{i=1}^{n+1} Vol(\mathcal{H}^h_i \cap \mathcal{S})}{Vol({\mathcal{S}})}. \tag{4.1}
\end{equation}
\end{definition}
It is clear, that $\sup_{\mathcal{S}}\delta(\mathcal{S}(h))$ provides an universal upper bound to any hyperball packing $\mathcal{H}^h$ in the
$n$-dimensional hyperbolic space
$\mathbb{H}^3$.   

{\it Due to Lemma \ref{lemma:3.6}, it suffices to consider regular truncated simplices and their corresponding densest hypersphere packings in the following.}

Let $\mathcal{S}$ be a regular $n$-simplex spanned by $n+1$ points (vertices) set $V=\{B_1,\dots,B_{n+1} \}$ in $\mathbb{H}^n$ where the vertices 
are outer point regarding the quadric $Q$.

The truncated regular tetrahedron $\mathcal{S}$ can be decomposed into $(n+1)!$ congruent simply truncated orthoschemes; one of them
$\mathcal{O}=Q_0Q_1Q_2P_0P_1P_2$ is illustrated in Fig.\ref{Fig.1}.a in $\HYP$ where $P_0$ is the centre of the "regular tetrahedron" $\mathcal{S}$,
$P_1$ is the centre of a hexagonal face of $\mathcal{S}$, $P_2$ is the midpoint of a ``common perpendicular" edge of this face,
$Q_0$ is the centre of an adjacent regular triangle face of $\mathcal{S}$, $Q_1$ is the midpoint of an appropriate edge of this face and
one of its endpoints is $Q_2$. 

Similarly, the higher-dimensional truncated regular simplex can be decomposed into characteristic truncated orthoschemes.
In $n$-dimensional hyperbolic space $\mathbb{H}^n$ $(n \ge 3)$ it can be seen that if $\mathcal{O}=Q_0Q_1Q_2 \dots Q_{n-1}$ $P_0P_1P_2 \dots P_{n-1}$ is a complete
orthoscheme  with degree $d=1$ (a simply frustum orthoscheme) where $B_1$ is a outer vertex of
$\overline{\mathbb{H}}^n$ then the points $Q_0,Q_1,Q_2,\dots,Q_{n-1}$ lie on the polar hyperplane $\beta_1$ of $B_1$ (see Fig.\ref{Fig.1} in $\HYP$).
The regular truncated simplices in $n$-dimensional hyperbolic space $\HYN$ can be constructed by orthoschemes of degree $1$ 
where one of their initial vertex (in the Fig.\ref{Fig.1}) $B_1$ is outer point regarding the quadric $Q$. 
These orthoschemes are characterized by their Coxeter-Schl\"afli graphs.
\subsection{Orthoschemes and the volume of a truncated regular tetrahedron}

{\it An orthoscheme $\mathcal{O}$ in $\mathbb{H}^n$ $n \geq 2$ in classical sense}
is a simplex bounded by $n+1$ hyperplanes $H_0,\dots,H_n$
such that $H_i \bot H_j, \  \text{for} \ j\ne i-1,i,i+1$ (\cite{K98, K91, Sz13-4}).
\begin{rmrk}
This definition is equivalent to the following: A simplex $\mathcal{O}$ in $\mathbb{H}^n$ is a
orthoscheme iff the $n+1$ vertices of $\mathcal{O}$ can be
labelled by $R_0,R_1,\dots,R_n$ in such a way that
$\text{span}(R_0,\dots,R_i) \perp \text{span}(R_i,\dots,R_n) \ \ \text{for} \ \ 0<i<n-1.$
\end{rmrk}

{\it The orthoschemes of degree} $d\in{0,1,2}$ in $\mathbb{H}^n$ are bounded by $n+d+1$ hyperplanes
$H_0,H_1,\dots,H_{n+d}$ such that $H_i \perp H_j$ for $j \ne i-1,~i,~i+1$, where, for $d=2$,
indices are taken modulo $n+3$.

Geometrically, complete orthoschemes of degree $m$ can be described as follows:
\begin{enumerate}
\item
For $d=0$, they coincide with the class of classical orthoschemes introduced by
{{Schl\"afli}}.
We denote the $(n+1)$-hyperface opposite to the vertex $R_i$
by $H_i$ $(0 \le i \le n)$. An orthoscheme $\mathcal{O}$ has $n$ dihedral angles which
are not right angles. Let $\alpha_{ij}$ denote the dihedral angle of $\mathcal{O}$
between the faces $H_i$ and $H_j$. Then we have
\begin{equation}
\alpha_{ij}=\frac{\pi}{2}, \ \ \text{if} \ \ 0 \le i < j -1 \le n. \notag
\end{equation}
The $n$ remaining dihedral angles $\alpha_{i,i+1}, \ (0 \le i \le n-1)$ are called the
essential angles of $\mathcal{O}$.
The initial and final vertices, $R_0$ and $R_n$ of the orthogonal edge-path
$R_iR_{i+1},~ i=0,\dots,n-1$, are called principal vertices of the orthoscheme.
\item
A complete orthoscheme of degree $d=1$ can be interpreted as an
orthoscheme with one outer principal vertex, say $R_n$, which is truncated by
its polar plane $pol(R_n)$ (see Fig.~4.~b). In this case the orthoscheme is called simply truncated with
ideal vertex $R_0$.
\item
A complete orthoscheme of degree $d=2$ can be interpreted as an
orthoscheme with two outer principal vertex, $R_0,~R_n$, which is truncated by
its polar hyperplanes $pol(R_0)$ and $pol(R_n)$. In this case the orthoscheme is called doubly
truncated. (In this case we distinguish two different type of the orthoschemes but we
shall not deeply go into the details (see \cite{K98, K91, Sz13-4}).
\end{enumerate}
The problem of determining $\sup_{\mathcal{S}}\delta(\mathcal{S})$ seems to be complicated in general, but we investigate this question in the $4$-dimensional hyperbolic space.
By the Lemma \ref{lemma:3.6}, it is sufficient to determine
the density upper bound of congruent saturated hyperball packings in for regular truncated simplices. 

First we have to determine the volumes of "$4$-dimensional simple frustum orthoschemes". 

\begin{figure}[h!]
\centering
\includegraphics[width=0.7\textwidth]{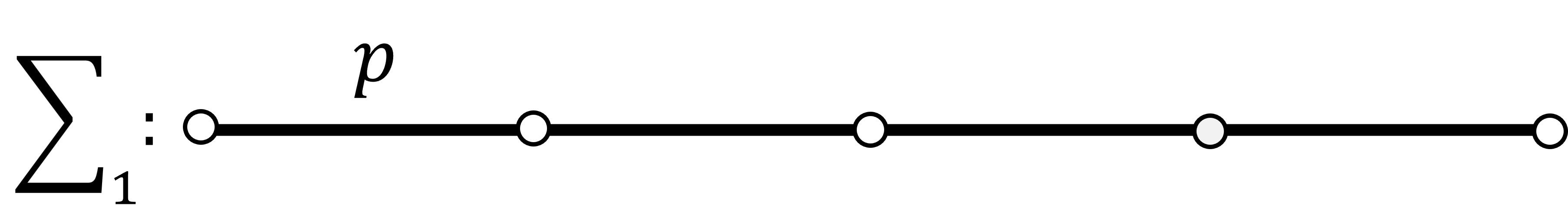}
\caption{Coxeter graph of the orthoscheme from dissected simplex $\mathcal{S}$. The possible values $p$ are provided in Lemma (\ref{Lemma:p_value}).}
\label{Coxeter_graph}
\end{figure}

The regular (truncated) $n$-simplex in $\mathbb{H}^n$ can be dissected into $(n+1)!$ congruent (truncated) orthoschemes. Thus, $\mathcal{S}(B_1,\dots,B_5)$ 
in $\mathbb{H}^4$ is dissected into $5!=120$ orthoschemes whose the Coxeter graph is given by Fig.~\ref{Coxeter_graph}.
The volume of an orthoscheme in even dimensions is computed by  Schl\"{a}fli reduction formula which is revisited and 
extended for complete orthoschemes of degree $1$ 
by R. Kellerhals in \cite{K91}.
\subsection{The Schl\"{a}fli function and reduction formula}
Following the Kellerhals'-Schl\"{a}fli's method (see \cite{K91}), we consider Schl\"{a}fli function, which plays a principal role in volume computations.
We start with Schl\"{a}fli function on spherical space. A spherical orthoscheme $R$ of dimension $n$ is represented by a linear elliptic scheme $\Sigma$ (Coxeter graph) of order $n+1$. The Schl\"{a}fli normalized volume function for $R$ is formulated by
\begin{equation*}
    f_n\left(\Sigma \right):=c_n\cdot \mathrm{vol}\left(R \right), \notag
\end{equation*}
with $\displaystyle c_n=\frac{2^{n+1}}{\mathrm{vol}\left(\mathbb{S}^n\right)}=\frac{2^n}{\pi^{\left(\frac{n+1}{2} \right)}}\Gamma\left( \frac{n+1}{2}\right)$, and $f_0:=1$.\\
As it stated in \cite{K91}, that function was generalized for complete orthoscheme in hyperbolic space. Let $\Sigma_{d}$ be Coxeter graph of a complete orthoscheme $R_d$ of degree $d, 0\leq d \leq 2$ in $\mathbb{H}^n$. 
The Schl\"{a}fli's function $F$ for complete orthoscheme 
$R_d \subset \mathbb{H}^n$ with graph $\Sigma_d$ is defined by
\begin{equation}
    F_n\left(\Sigma_d\right):=i^{n}\cdot c_n\cdot\mathrm{vol}_n\left(R_d \right), \tag{4.2}
\end{equation}
where $F_0:=1$, $i^2=-1$, and 
\begin{equation*}
    c_n=\frac{2^n}{\pi^{\left(\frac{n+1}{2} \right)}}\Gamma\left( \frac{n+1}{2}\right). \tag{4.3}
\end{equation*}
Furthermore, the following Schl\"{a}fli reduction formula is the core in our computation.
\begin{theorem}[Schl\"{a}fli reduction formula \cite{K91}]\label{thm:reduction_schlafli}
Denote by $R_d \subset \mathbb{H}^{2n}$, $0\leq d \leq 2$, $n \geq 1$, a $2n$-dimensional orthoscheme of degree $d$ with scheme $\Sigma_{d}$. Then
\begin{equation}
F_{2n}\left( \Sigma_{d} \right)=\sum_{k=1}^{n}\frac{(-1)^k}{k+1}\begin{pmatrix}2k \\ k \end{pmatrix}\sum_{\sigma}f_{2n-(2k+1)}\left(\sigma \right), \tag{4.4}
\end{equation}
where $\sigma$ runs through all elliptic subschemes of order $2(n-k)$ of $\Sigma_d$ all of whose components are even order.
\end{theorem}
According to Theorem~\ref{thm:reduction_schlafli}, and our orthoscheme, we set $n=2$ (having $2n=4$ dimensional) and $d=1$ (simply truncated orthoscheme). 
\begin{equation*}
F_4\left(\Sigma_1\right)=\sum_{\sigma_4}f_3(\sigma_4)-\frac{1}{2}\begin{pmatrix} 2 \\ 1 \end{pmatrix}\sum_{\sigma_2}f_1(\sigma_2)+\frac{1}{3}\begin{pmatrix} 4 \\ 
2 \end{pmatrix}\sum_{\sigma_0}f_{-1}(\sigma_0),
\end{equation*}
where $\sigma_4,~\sigma_2,~\sigma_0$ are subschemes of order $4,~2,~0$ respectively.\\
We compute, $F_4\left(\Sigma_1\right)=\frac{2}{15}-\frac{2}{3p}$. On the other hand $F_4\left(\Sigma_1\right)=(-1)^{4} c_4Vol_{4}(R_1)$, 
where $c_4=\frac{12}{\pi^2}$. Thus, we obtain the following 
\begin{lemma} 
The volume $Vol_4(R_1)$ of the orthoscheme $R_1$ related to the Coxeter graph $\Sigma_1$ (see Fig.~\ref{Coxeter_graph}) is
\begin{equation}
Vol_4(R_1)=\frac{F_4(\Sigma_1)}{c_4}=\frac{\pi^2}{12}\left(\frac{2}{15}-\frac{2}{3p} \right). \tag{4.5}
\end{equation}
\end{lemma} 
\begin{figure}[h!]
	\centering
	\includegraphics[width=0.75\textwidth]{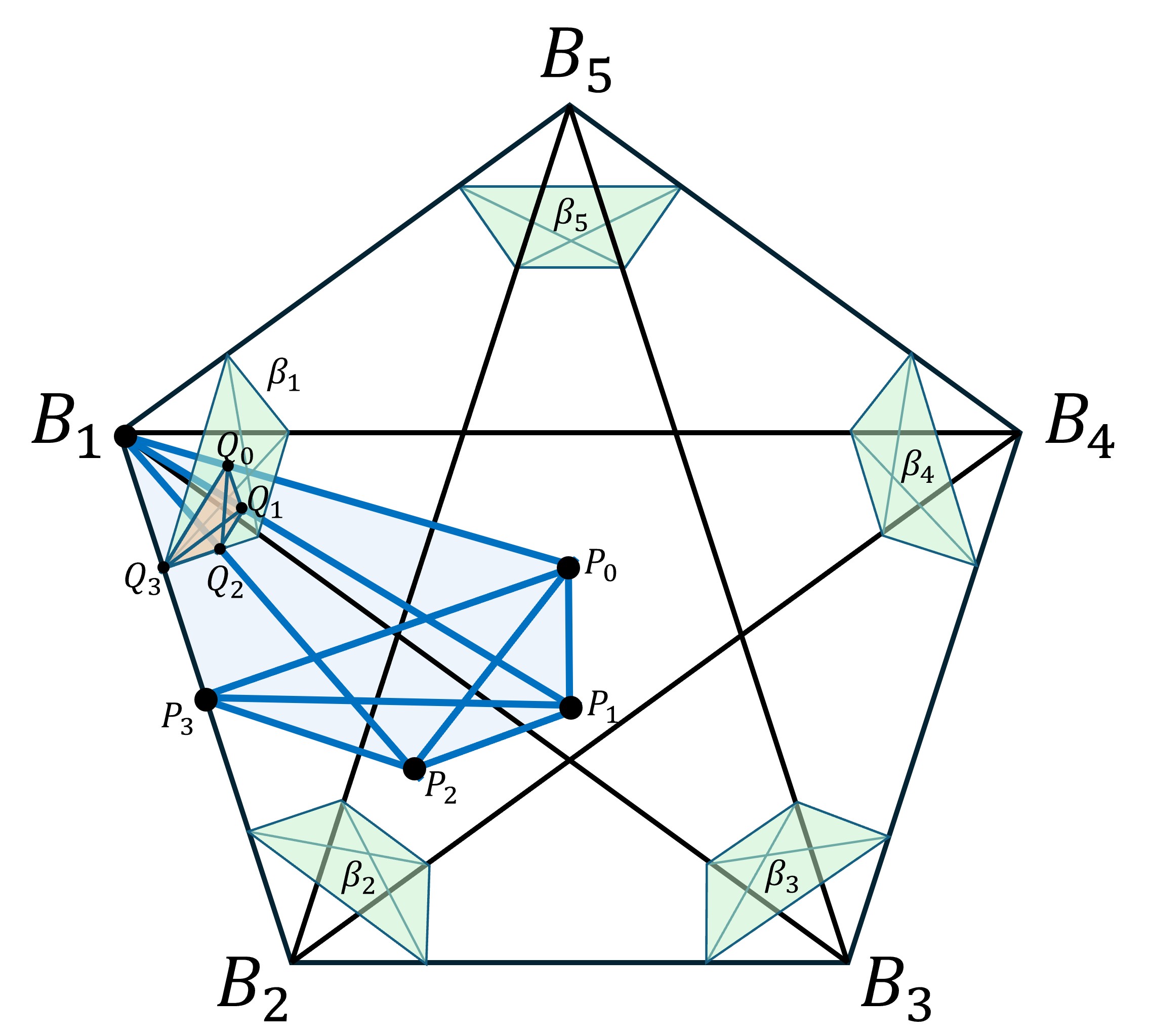}
	\caption{The structure of orthoscheme obtained from the regular simplex in $\mathbb{H}^4$}
	\label{Fig:dissection_method_2}
\end{figure}
We also need the volumes of the hyperball pieces lying in the truncated simplex, for which we use the apparatus of the Beltrami-Cayley-Klein model 
(see e.g., \cite{Sz17-3}).

By considering regular $4$-simplex $\mathcal{S}(B_1,B_2,B_3,B_4,B_5)$ in $\mathbb{H}^4$ whose vertices $B_1[\textbf{b}_1],~B_2[\textbf{b}_2],~B_3[\textbf{b}_3],~ B_4[\textbf{b}_4]$, and $B_5[\textbf{b}_5]$. Here, we choose the points $\mathbf{b}_i,~i=1, \dots, 5$ in the Beltrami-Cayley-Klein model, such that the center of $\mathcal{S}(B_1,B_2,B_3,B_4,B_5)$, $P_{0}$, is exactly at $(1,0,0,0,0)$:\\
$\mathbf{b}_1=\left(1,0,0,0,s \right)$, $\mathbf{b}_2=\left(1,0,\frac{\sqrt{10}}{4}s,\frac{\sqrt{5}}{4}s,-\frac{1}{4}s \right)$, $\mathbf{b}_3=\left(1,0,-\frac{\sqrt{10}}{4}s,-\frac{\sqrt{5}}{4}s,-\frac{1}{4}s \right)$,  \\
$\mathbf{b}_4=\left(1,-\frac{\sqrt{10}}{4}s,0,-\frac{\sqrt{5}}{4}s,-\frac{1}{4}s \right)$, $ \mathbf{b}_5=\left(1,\frac{\sqrt{10}}{4}s,0,-\frac{\sqrt{5}}{4}s,-\frac{1}{4}s\right)$,  where $s\geq 1$. In particular, if $s=1$, then $\mathcal{S}$ is an ideal simplex whose all vertices $B_i$ lie exactly at the infinity.\\
By letting $s>1$, $\mathcal{S}$'s vertices are ultra-ideal (lying outer of the BCK model). As a consequence, we can construct a polar hyperplanes of these ultra-ideal vertices, denoted by $\beta_i:=pol(B_i)$. Moreover, one could obtain the compact polyhedron by truncating the simplex $\mathcal{S}$ with hyperplanes $\beta_i$ associated with its vertices, (see Fig.~\ref{Coxeter_graph}.a. and Fig.~\ref{Fig:dissection_method_2}). 


We dissect the truncated simplex $\mathcal{S}(B_1,B_2,B_3,B_4,B_5)$ 
into characteristic orthoschemes by usual way.

Note that in our coordinate system related to the vectors $\mathbf{b}_1 \dots \mathbf{b}_5$, it is simply determined that the center coordinate of the regular simplex $\mathcal{S}$, denoted by $P_0$, is given by $P_0[\mathbf{p_0}]$, where $\displaystyle\mathbf{p_0}=\frac{1}{5}\sum_{i=1}^{5}\mathbf{b}_i$. It also holds for other centers of lower-dimensional faces of $\mathcal{S}$ denoted by $P_1$, $P_2$, $P_3$ (the centers of faces $B_1B_2B_3B_4$, $B_1B_2B_3$, $B_1B_2$, respectively). In our model, these centers have coordinates $P_1[\mathbf{p_1}], P_2[\mathbf{p_2}], P_3[\mathbf{p_3}]$, where 
\begin{equation}
\begin{gathered}
\mathbf{p_1}=\frac{1}{4}\sum_{i=1}^{4}\mathbf{b}_i, ~ \mathbf{p_2}=\frac{1}{3}\sum_{i=1}^{3}\mathbf{b}_i, ~ \mathbf{p_3}=\frac{1}{2}\sum_{i=1}^{2}\mathbf{b}_i, \\
\mathbf{p}_0=\begin{pmatrix}1 & 0 &0 &0 &0\end{pmatrix},~ \mathbf{p_1}=\begin{pmatrix}1 &-\frac{s\sqrt{10}}{16} &0& \frac{s\sqrt{5}}{16} & \frac{s}{16} \end{pmatrix},\\ \mathbf{p}_2=\begin{pmatrix}1 &0 &0 & \frac{s\sqrt{5}}{16} & \frac{s}{6} \end{pmatrix},~ \mathbf{p}_3=\begin{pmatrix}1 & 0 & \frac{s \sqrt{10}}{8}& \frac{s\sqrt{5}}{8} & \frac{3s}{8} \end{pmatrix}, 
\end{gathered} \notag
\end{equation}

Moreover, the associated (normalized) forms of the faces of the orthoscheme $P_0P_1P_2P_3B_1$ are given by
\begin{align*}
\Bb^{1}=\begin{bmatrix} \frac{s}{\sqrt{16-s^2}} \\ \frac{\sqrt{10}}{\sqrt{16-s^2}}\\0\\\frac{\sqrt{5}}{\sqrt{16-s^2}}\\\frac{1}{\sqrt{16-s^2}}\end{bmatrix},~~\Bb^{2}=\begin{bmatrix}
0 \\ 1 \\ 0 \\ 0 \\ 0
\end{bmatrix},~~\Bb^{3}=\begin{bmatrix}
0 \\ -\frac{1}{2} \\ \frac{1}{2} \\ \frac{\sqrt{2}}{2} \\ 0
\end{bmatrix},~~\Bb^{4}=\begin{bmatrix}
0 \\ 0 \\ -1 \\ 0 \\ 0
\end{bmatrix},~~\Bb^{5}=\begin{bmatrix}
0 \\ 0 \\ \frac{1}{2} \\ \frac{\sqrt{2}}{4} \\ -\frac{\sqrt{10}}{4}
\end{bmatrix}.
\end{align*}
The corresponding Coxeter-Schl\"{a}fli matrix formed by Grammanian of these forms is given by
\begin{equation}\label{Coxeter-Schlafli_Matrix_Orthoscheme}
G:=\left(\langle \Bb^{i},\Bb^{j} \rangle \right)=\begin{pmatrix}
1 & -\frac{\sqrt{10}}{\sqrt{16-s^2}} & 0 & 0 & 0 \\
-\frac{\sqrt{10}}{\sqrt{16-s^2}} & 1 & -\frac{1}{2} & 0 & 0 \\
0 & -\frac{1}{2} & 1 & -\frac{1}{2} & 0 \\
0 & 0 & -\frac{1}{2} & 1 & -\frac{1}{2} \\
0 & 0 & 0 & -\frac{1}{2} & 1
\end{pmatrix} \tag{4.6}
\end{equation}
Let $\frac{\pi}{p}$ be the angle between two hyperplanes of forms $\Bb^1$ and $\Bb^2$. So that the $-\cos{\left(\frac{\pi}{p}\right)}$ is equal to the element at the second row and first column $G$, i.e.,
\begin{equation}\label{Eq:4.7}
\cos{\left(\frac{\pi}{p}\right)}=\frac{\sqrt{10}}{\sqrt{16-s^2}}\Rightarrow s=\frac{\sqrt{16\cos^{2}{\left( \frac{\pi}{p}\right)}-10}}{\cos^2{\left( \frac{\pi}{p} \right)}}, \frac{\pi}{p}=\arccos{\left(\sqrt{\frac{10}{16-s^2}}\right)} . \tag{4.7}
\end{equation}
 In particular, if $s=1$ (the principal vertex $B_1$ is an ideal vertex, as well as other vertices of the simplex $\mathcal{S}$), then we have $p=\frac{\pi}{\arccos{\left( \sqrt{\frac{2}{3}} \right)}}$. Meanwhile, if $p>\frac{\pi}{\arccos{\left( \sqrt{\frac{2}{3}} \right)}}$, then $B_1$ (and all vertices of the simplex are ultra ideal).\\
 Thus, the polar hyperplane $\beta_1:=pol(B_1)$ truncates this orthoscheme (and so the simplex $\mathcal{S}$).\\
Meanwhile, the simplex is truncated by other vertices polar hyperplanes $\beta_i$.\\
For each faces $\beta_i$ we construct congruent hyperballs of height $h>0$. We maximize $h$ as a half of the distance of any two polar hyperplanes, i.e., $h=\frac{1}{2}d(\beta_i,\beta_j)$, $i\neq j$.\\
Therefore, we also require $d(\beta_i,\beta_j)>0$, that is $\beta_i$, $\beta_j$, intersect at ultra ideal points,
\begin{equation}\label{Eq:require_h}
\left|\frac{\langle \mathbf{b}_i, \mathbf{b}_j \rangle}{\sqrt{\langle \mathbf{b}_i, \mathbf{b}_i \rangle \langle \mathbf{b}_j, \mathbf{b}_j \rangle}}\right| \geq 1. \tag{4.8}
\end{equation}

Without loss of generality, we choose $\mathbf{b}_1$ and $\mathbf{b}_2$, and substitute them together with (\ref{Eq:4.7}) into (\ref{Eq:4.9}) we conclude that $\displaystyle \cos \left(\frac{\pi}{p}\right)< \sqrt{\frac{3}{4}}$. Hence, we obtain the following
\begin{lemma}\label{Lemma:p_value}
In regular $4$-simplex $\mathcal{S}(B_1,B_2,B_3,B_4,B_5)$, with outer vertices is determined by the parameter $p$ where it has to satisfy the following conditions
\begin{equation}
\frac{\pi}{\arccos{\left( \sqrt{\frac{2}{3}}\right)}} < p < \frac{\pi}{\arccos{\left( \sqrt{\frac{3}{4}}\right)}}. \notag
\end{equation}
($5.1042\dots  < p < 6$).
\end{lemma}
\begin{rmrk}
The dihedral angles of the simplex is $\frac{2\pi}{p}$, that is exactly twice of the essential angle of the orthoscheme whose the form $\frac{\pi}{p}$.
Let $\alpha$ be the dihedral angles of the original simplex $\mathcal{S}(B_1B_2B_3B_4B_5)$. The corresponding Coxeter-Schl\"{a}fli matrix of $\mathcal{S}$ is given by
\begin{equation*}
G_{\mathcal{S}}:=\begin{pmatrix}
1 & \frac{s^2+4}{s^2-16} & \frac{s^2+4}{s^2-16} & \frac{s^2+4}{s^2-16} & \frac{s^2+4}{s^2-16}\\
\frac{s^2+4}{s^2-16} & 1  & \frac{s^2+4}{s^2-16} & \frac{s^2+4}{s^2-16} & \frac{s^2+4}{s^2-16}\\
\frac{s^2+4}{s^2-16} & \frac{s^2+4}{s^2-16} & 1  & \frac{s^2+4}{s^2-16} & \frac{s^2+4}{s^2-16} \\
  \frac{s^2+4}{s^2-16} & \frac{s^2+4}{s^2-16} & \frac{s^2+4}{s^2-16} & 1 & \frac{s^2+4}{s^2-16}  \\
 \frac{s^2+4}{s^2-16} & \frac{s^2+4}{s^2-16} & \frac{s^2+4}{s^2-16} & \frac{s^2+4}{s^2-16} & 1
\end{pmatrix}.
\end{equation*}

Thus, $\frac{s^2+4}{s^2-16}=-\cos{\alpha}$. By expressing $s$ in $p$ based on relation (\ref{Eq:4.7}), we obtain $\cos{\alpha}=2\cos^2{\frac{\pi}{p}}-1=\cos{\frac{2\pi}{p}}$.
\end{rmrk}
\subsection{Volume of hyperball part}
The volume of hyperball part with base $\mathcal{A}_3$ on a polar hyperplane, and height $h$ in $\mathbb{H}^4$ is given by formula (\ref{Eq:2.2}) {here, we set~$k=\sqrt{\frac{-1}{K}}=1$.

{\bf The maximum possible height} $h(p)$ of the hyperballs belonging to the simplex $\mathcal{S}(B_1,B_2,B_3,B_4,B_5)$ with base planes $\beta_i$ are exactly the half of the distance between any pair of polar hyperplane, precisely
\begin{equation}\label{Eq:4.9}
\begin{aligned}
h(p)&=\frac{1}{2}d(\beta_i,\beta_j)=\frac{1}{2}{\mathrm{arccosh}}{\left(-\frac{\langle \boldsymbol{b}_i, \boldsymbol{b}_j \rangle}{\sqrt{\langle \boldsymbol{b}_i, \boldsymbol{b}_i 
\rangle \langle \boldsymbol{b}_j, \boldsymbol{b}_j \rangle}}\right)}\\&=\frac{1}{2}{\mathrm{arccosh}}{\left(\frac{s^2+4}{4s^2-4}\right)}=\frac{1}{2}{\mathrm{arccosh}}
{\left(\frac{\cos{\left(\frac{2\pi}{p} \right)}}{3\cos{\left(\frac{2\pi}{p}\right)}-1}\right)}.
\end{aligned} \tag{4.9}
\end{equation}
The base $\mathcal{A}_3$ lies on $\beta_i$, can be obtained by considering the lines connecting the vertices of orthoschemes $P_0P_1P_2P_3$ to the ultra ideal point $B_1$. These lines intersect the polar hyperplane $\beta_1$ at namely $Q_0$, $Q_1$, $Q_2$, $Q_3$ respectively. These are $Q_0:=\overline{B_1P_0}\cap\beta_1$, $Q_1:=\overline{B_1P_1}\cap\beta_1$, $Q_2:=\overline{B_1P_2}\cap\beta_1$, $Q_3:=\overline{B_1P_3}\cap\beta_1$.
Based on some direct computations we obtain the Coxeter Schl\"{a}fli matrix of polyhedron $\mathcal{A}_3$ in $\mathbb{H}^3$ as follows
\begin{equation}
G_{\mathcal{A}_3}=\begin{pmatrix}
1& -\cos{\left(\frac{\pi}{p}\right)}& 0 & 0 \\ -\cos{\left(\frac{\pi}{p}\right)} & 1 & -\frac{1}{2} & 0 \\
0 & -\frac{1}{2} & 1 & -\frac{1}{2} \\
0 & 0 & -\frac{1}{2} & 1
\end{pmatrix}. \tag{4.10}
\end{equation}
Thus, the polyhedron $\mathcal{A}_3$ is a 3-orthoscheme. 
The volume of $\mathcal{A}_3$ can be directly computed by the Kellerhals formula \cite{K98}:
\begin{theorem}\label{Vol_Kellerhals} The volume of a three-dimensional hyperbolic
complete ortho\-scheme (e.g., without Lambert cube cases, where $a_0, a_3$ intersect each other) $\mathcal{R}$
is expressed with the essential angles $\alpha^{01},\alpha^{12},\alpha^{23}, (0 \le \alpha^{ij} \le \frac{\pi}{2})$ in the following:
\begin{align*}
&Vol_3(\mathcal{R})=\frac{1}{4} \{ \mathcal{L}(\alpha^{01}+\theta)-
\mathcal{L}(\alpha^{01}-\theta)+\mathcal{L}(\frac{\pi}{2}+\alpha^{12}-\theta)+ \notag \\
&+\mathcal{L}(\frac{\pi}{2}-\alpha^{12}-\theta)+\mathcal{L}(\alpha^{23}+\theta)-
\mathcal{L}(\alpha^{23}-\theta)+2\mathcal{L}(\frac{\pi}{2}-\theta) \}, 
\end{align*}
where $\theta \in [0,\frac{\pi}{2})$ is defined by the following formula:
$$
\tan(\theta)=\frac{\sqrt{ \cos^2{\alpha^{12}}-\sin^2{\alpha^{01}} \sin^2{\alpha^{23}
}}} {\cos{\alpha^{01}}\cos{\alpha^{23}}}
$$
and, where $\mathcal{L}(x):=-\int\limits_0^x \log \vert {2\sin{t}} \vert dt$ \ denotes the
Lobachevsky function.
\end{theorem}
As a consequence, applying this for $Q_0Q_1Q_2Q_3$ orthoscheme we obtain the explicit expression
\begin{equation}\label{Vol_A3}
\begin{aligned}
Vol_3(\mathcal{A}_3)=&\frac{1}{4} \{ \mathcal{L}(\frac{\pi}{p}+\theta)-
\mathcal{L}(\frac{\pi}{p}-\theta)+\mathcal{L}(\frac{\pi}{2}+\frac{\pi}{3}-\theta)+ \notag \\
&+\mathcal{L}(\frac{\pi}{2}-\frac{\pi}{3}-\theta)+\mathcal{L}(\frac{\pi}{3}+\theta)-
\mathcal{L}(\frac{\pi}{3}-\theta)+2\mathcal{L}(\frac{\pi}{2}-\theta) \},\\
&\text{where, } \theta:=\theta(p)=\arctan{\left(\frac{\sqrt{1-3\sin^2{\left( \frac{\pi}{p}\right)}}}{\cos{\left(\frac{\pi}{p} \right)}} \right)}.
\end{aligned} \tag{4.11}
\end{equation}
In essence, the volume of base $\mathcal{A}_3$ is completely depend on the parameter $p$. It is then reasonable to write $Vol_3(\mathcal{A}_3):=Vol_3\left(\mathcal{A}_3(p)\right)$.
\subsection{The local density of hyperball packings related to simplex $\mathcal{S}(B_1,B_2,B_3,B_4,B_5)$ }
The density of hyperball packing in a regular truncated simplex $\mathcal{S}(B_1,B_2,B_3,B_4,B_5)$ } can be defined as a comparison between the volume sum of hyperball parts lying in the simplex and the volume of the truncated simplex. It is clear, that it is equal to the ratio of the volume of the hybeball part contained by the $R_1$ orthoscheme to the volume of the $R_1$ orthoscheme.
\begin{definition}\label{def:density}
 The local density $\delta(p)$ of the congruent hyperball packings related to the regular truncated simplices is 
 \begin{equation}
\delta{(p)}=\frac{Vol_4(\mathcal{H}^{h^+(p)}_4(\mathcal{A}_3(p)))}{Vol_4\left(R_1(p) \right)}, \notag
\end{equation}
\begin{center}
where $\displaystyle \frac{\pi}{\arccos{\left( \sqrt{\frac{2}{3}}\right)}} < p < \frac{\pi}{\arccos{\left( \sqrt{\frac{3}{4}}\right)}}$.
\end{center}
\end{definition}
Using Definition \ref{def:density} and formulas (\ref{Eq:4.9}), (\ref{Vol_A3})  and Lemma \ref{Lemma:p_value} we obtain the following 
\begin{theorem}\label{thm:local_density}
The locally optimal congruent hyperball packing density related to the regular truncated simplices for a given possible parameter $p$ is: 
\begin{equation}
\delta{(p)}=\frac{\frac{1}{8}Vol_3(\mathcal{A}_3(p))\cdot k \cdot \left[\frac{2}{3}\sinh{\left(\frac{3h(p)}{k}\right)}+
6\sinh{\left(\frac{h(p)}{k} \right)}\right]}{\frac{\pi^2}{12}\left(\frac{2}{15}-\frac{2}{3p}\right)} \notag
\end{equation}
where $\displaystyle h(p)=\frac{1}{2}{\mathrm{arccosh}}{\left(\frac{\cos{\left(\frac{2\pi}{p} \right)}}{\left\lvert 3\cos{\left(\frac{2\pi}{p}\right)}-1\right\rvert}\right)}$ and $\frac{\pi}{\arccos{\left( \sqrt{\frac{2}{3}}\right)}} < p < \frac{\pi}{\arccos{\left( \sqrt{\frac{3}{4}}\right)}}$. Here we set the constant curvature $k=\sqrt{\frac{-1}{K}}=1$. \quad \quad $\square$
\end{theorem} 

To determine the densest possible hyperball arrangement for regular truncated 4-dimensional simplices, we need to examine the density function defined in Theorem \ref{thm:local_density}, which depends on the variable $p$. 

\begin{figure}[h!]
\centering
\includegraphics[width=\textwidth]{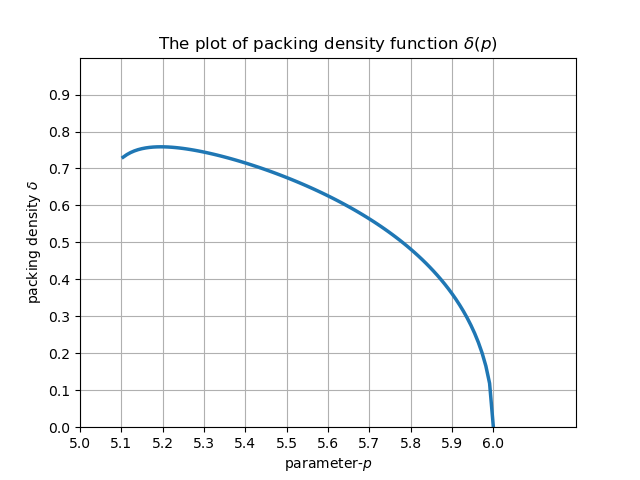}
\caption{The plot of density function $\delta(p)$}
\label{Fig:4.density}
\end{figure}
Finally, we obtain the plot after careful analysis of the smooth
density function (cf. Fig.~\ref{Fig:4.density}) and we obtain the following
\begin{theorem}\label{thm:4.10:density_function}
The density function $\delta(\mathcal{S}(p))$, $p\in $
attains its maximum at $p^{opt} \approx 5.19550$, and $\delta(\mathcal{S}(p))$
is strictly increasing in the interval $(\frac{\pi}{\arccos{\sqrt{\frac{2}{3}}}},p^{opt})$, and strictly decreasing in $(p^{opt},\frac{\pi}{\arccos{\sqrt{\frac{3}{4}}}})$. Moreover, the optimal density
$\delta^{opt}(\mathcal{S}(p^{opt})) \approx 0.7586482$ (see Fig.~\ref{Fig:4.density}).
\end{theorem}
\begin{rmrk}
\begin{enumerate}
\item In our case $\lim_{p\rightarrow \frac{\pi}{\arccos{\left( \sqrt{\frac{2}{3}}\right)}}}(\delta(\mathcal{S}(p)))$ is equal to the conjectured B\"{o}r\"{o}czky-Florian
upper bound of the ball and horoball packings in $\mathbb{H}^4$ \cite{B78,K98}.
\item $\delta^{opt}(\mathcal{S}(p^{opt})) \approx 0.7586482$ is larger than
the B\"{o}r\"{o}czky-Florian type conjectured upper bound
$\delta_{BF} \approx 0.73....$; but these hyperball packing configurations
are only locally optimal and cannot be extended to the entire hyperbolic
space $\mathbb{H}^4$.
\end{enumerate}
\end{rmrk}
We obtain the next theorem as the direct consequence of the previous statements:
\begin{theorem}
The density upper bound of the saturated congruent hyperball packings related to the
corresponding truncated simplex  cells is realized in a regular truncated tetrahedra belonging to parameter $p^{opt} \approx 5.19550$ with density $\approx 0.7586482$.
\end{theorem}
It follows immediately from the analysis of the density function $\delta(p)$, but the graph (see Fig.~4) of the function also clearly shows that it is not monotonic.
We get from the Theorem \ref{thm:4.10:density_function} directly the denial of the A.~Przeworski's conjecture \cite{P13}:
\begin{cor}
The density function $\delta(\mathcal{S}(p))$, is not an increasing function of $h(p)$ (the height of hyperballs).
\end{cor}
%


%
\noindent
\footnotesize{Budapest University of Technology and Economics, Institute of Mathematics, \\
Department of Geometry, \\
H-1521 Budapest, Hungary. \\
E-mail:~szirmai@math.bme.hu \\
http://www.math.bme.hu/ $^\sim$szirmai}

\end{document}